\documentclass[11pt]{amsart}
\usepackage{amsmath,amsfonts,amssymb}
\usepackage{stmaryrd}
\usepackage[colorlinks=true, linkcolor=blue,
  citecolor=blue, urlcolor=blue]{hyperref}

\def\newthm#1#2{\newtheorem{#1}[dummy]{#2}%
  \expandafter\def\csname#2\endcsname##1{\hyperref[#1:##1]{#2~\ref*{#1:##1}}}}
\newthm{thm}{Theorem}
\newthm{lemma}{Lemma}
\newthm{prop}{Proposition}
\newthm{cor}{Corollary}
\newthm{conj}{Conjecture}
\newthm{cons}{Consequence}
\newthm{fact}{Fact}
\newthm{obs}{Observation}
\newthm{guess}{Guess}
\theoremstyle{definition}
\newthm{defn}{Definition}
\newthm{remark}{Remark}
\newthm{example}{Example}
\newthm{question}{Question}
\newthm{notation}{Notation}

\newcommand{\Section}[1]{\hyperref[sec:#1]{Section~\ref*{sec:#1}}}
\newcommand{\Table}[1]{\hyperref[tab:#1]{Table~\ref*{tab:#1}}}
\newcommand{\eqn}[1]{\hyperref[eqn:#1]{(\ref*{eqn:#1})}}

\DeclareMathOperator{\Gr}{\textrm{Gr}}

\DeclareMathOperator{\LG}{\textrm{LG}}

\DeclareMathOperator{\OG}{\textrm{OG}}

\DeclareMathOperator{\QH}{\textrm{QH}}
\DeclareMathOperator{\QK}{\textrm{QK}}

\DeclareMathOperator{\codim}{\textrm{codim}}

\DeclareMathOperator{\dist}{dist}
\def\poly{{\mathrm{poly}}}

\usepackage{bm}

\newcommand{\bP}{{\mathbb P}}

\newcommand{\C}{{\mathbb C}}

\newcommand{\Z}{{\mathbb Z}}

\newcommand{\cO}{{\mathcal O}}

\newcommand{\euler}[1]{\chi_{_{#1}}}
\newcommand{\pt}{\textrm{pt}}

\newcommand{\al}{{\alpha}}

\newcommand{\ga}{{\gamma}}

\newcommand{\ev}{\textrm{ev}}

\newcommand{\wh}{\widehat}
\newcommand{\wb}{\overline}
\newcommand{\ov}{\overline}

\newcommand{\ignore}[1]{}

\newcommand{\Mb}{\wb{\mathcal M}}

\begin{document}

\title{Euler characteristics of cominuscule quantum $K$-theory}

\date{November 18, 2016}

\author{Anders~S.~Buch}
\address{Department of Mathematics, Rutgers University, 110
  Frelinghuysen Road, Piscataway, NJ 08854, USA}
\email{asbuch@math.rutgers.edu}

\author{Sjuvon Chung}
\address{Department of Mathematics, Rutgers University, 110
  Frelinghuysen Road, Piscataway, NJ 08854, USA}
\email{sjuvon.chung@rutgers.edu}


\thanks{Both authors were supported in part by NSF grants
  DMS-1503662.}

\begin{abstract}
  We prove an identity relating the product of two opposite Schubert
  varieties in the (equivariant) quantum $K$-theory ring of a
  cominuscule flag variety to the minimal degree of a rational curve
  connecting the Schubert varieties.  We deduce that the sum of the
  structure constants associated to any product of Schubert classes is
  equal to one.  Equivalently, the sheaf Euler characteristic map
  extends to a ring homomorphism defined on the quantum $K$-theory
  ring.
\end{abstract}

\maketitle

\section{Introduction}

Let $X = G/P$ be a flag variety defined by a semisimple complex Lie
group $G$ and a parabolic subgroup $P$.  The (small) equivariant
quantum $K$-theory ring $\QK_T(X)$ of Givental \cite{givental:wdvv} is
a formal deformation of the Grothendieck ring $K_T(X)$ of
$T$-equivariant algebraic vector bundles on $X$, where $T$ is a
maximal torus in $G$.  This ring encodes geometric information about
families of rational curves meeting triples of general Schubert
varieties in $X$, including the arithmetic genera of such families.
In this note we prove three results that present aspects of this
information in a concrete form when $X$ is a cominuscule flag variety,
that is, a Grassmannian $\Gr(m,n)$ of type A, a Lagrangian
Grassmannian $\LG(n,2n)$, a maximal orthogonal Grassmannian
$\OG(n,2n)$, a quadric hypersurface $Q^n$, or one of two exceptional
spaces known as the Cayley plane $E_6/P_6$ and the Freudenthal variety
$E_7/P_7$.

The ring $K_T(X)$ has a basis of Schubert structure sheaves $\cO^w =
[\cO_{X^w}]$ over the ring $\Gamma = K_T(\pt)$ of virtual
representations of $T$.  The quantum $K$-theory ring $\QK_T(X)$
consists of all formal power series with coefficients in $K_T(X)$.
The product of two Schubert classes in this ring has the form
\begin{equation}\label{eqn:strconst}
  \cO^u \star \cO^v = \sum_{w,d \geq 0} N^{w,d}_{u,v}\, q^d\, \cO^w \,,
\end{equation}
where the sum is over all Schubert classes $\cO^w$ and effective
degrees $d \in H_2(X;\Z)$.  Givental defined the structure constants
$N^{w,d}_{u,v} \in \Gamma$ as polynomial expressions in the
$K$-theoretic Gromov-Witten invariants of $X$ and proved that the
resulting product is associative \cite{givental:wdvv}.  The structure
of the ring $\QK_T(X)$ has been studied in the cominuscule case in a
series of papers by Chaput, Mihalcea, Perrin, and the first author
\cite{buch.mihalcea:quantum, chaput.perrin:rationality,
  buch.chaput.ea:finiteness, buch.chaput.ea:projected,
  buch.chaput.ea:chevalley}.  In particular, it has been proved that
only finitely many of the coefficients $N^{w,d}_{u,v}$ are non-zero
\cite{buch.chaput.ea:finiteness, buch.chaput.ea:rational}.  The
quantum $K$-theory of Grassmannians of type A has been related to
integrable systems by Gorbounov and Korff
\cite{gorbounov.korff:quantum}.  Conjectures for the ring structure of
$\QK_T(X)$ have also been given by Lenart and Maeno
\cite{lenart.maeno:quantum} and Lenart and Postnikov
\cite{lenart.postnikov:affine} when $X = G/B$ is defined by a Borel
subgroup of $G$.  In general the structure constants $N^{w,d}_{u,v}$
are conjectured to satisfy Griffeth-Ram positivity
\cite{griffeth.ram:affine}, that is, up to a sign these constants are
polynomials with non-negative coefficients in the classes
$[\C_{-\al}]-1 \in \Gamma$, where $\C_{-\al}$ is any one-dimensional
representation of $T$ defined by a negative root (see e.g.\
\cite{buch.chaput.ea:chevalley}).  This conjecture has been proved for
the structure constants $N^{w,0}_{u,v}$ of the equivariant $K$-theory
ring $K_T(X)$ by Anderson, Griffeth, and Miller
\cite{anderson.griffeth.ea:positivity}, and for the equivariant
quantum cohomology ring $\QH_T(X)$ by Mihalcea
\cite{mihalcea:positivity}.

Assume now that $X$ is a cominuscule flag variety.  Our work started
with the experimental observation that the sum of the structure
constants defining any product $\cO^u \star \cO^v$ of Schubert classes
in $\QK_T(X)$ is equal to 1.  This is our first result.

\begin{thm}\label{thm:sumcoef}
  For fixed $u, v$ we have $\,\displaystyle\sum_{w,d\geq 0}
  N^{w,d}_{u,v} = 1\,$ in $\Gamma$.
\end{thm}

Let $\euler{X} : K_T(X) \to \Gamma$ be the sheaf Euler characteristic
map, defined as the equivariant pushforward along the structure
morphism $X \to \{\pt\}$.  Equivalently, $\euler{X}$ is the unique
$\Gamma$-linear map defined by $\euler{X}(\cO^w) = 1$ for all $w$.
While this map is not a ring homomorphism unless $X$ is a single
point, \Theorem{sumcoef} is equivalent to the following statement.

\begin{thm}\label{thm:ringhom}
  Let $\QK_T^\poly(X) \subset \QK_T(X)$ be the subring of all finite
  power series.  There exists a unique ring homomorphism $\widehat\chi
  : \QK_T^\poly(X) \to \Gamma$ defined by $\widehat\chi(q)=1$ and
  $\widehat\chi(\cO^w) = 1$ for all $w$.
\end{thm}

Given two opposite Schubert varieties $X^u$ and $X_v$ in the
cominuscule flag variety $X$, let $\dist(X^u,X_v)$ denote the minimal
degree of a rational curve connecting these subvarieties.  This degree
is the smallest power of the deformation parameter $q$ that occurs in
the product $\cO^u \star \cO_v$, where $\cO_v = [\cO_{X_v}]$.  Let
$\chi : \QK_T(X) \to \Gamma\llbracket q \rrbracket$ denote the
$\Gamma\llbracket q \rrbracket$-linear extension of the sheaf Euler
characteristic map, defined by $\chi(\cO^w) = 1$ for all $w$.  Both of
the above theorems are consequences of the following identity.

\begin{thm}\label{thm:dist}
  We have $\chi(\cO^u \star \cO_v) = q^{\dist(X^u,X_v)}$.
\end{thm}

The proof of \Theorem{dist} is based on a construction of the ring
$\QK_T(X)$ using projected Gromov-Witten varieties
\cite{buch.chaput.ea:chevalley}, together with a relation between such
varieties and $K$-theoretic Gromov-Witten invariants
\cite{knutson.lam.ea:positroid, buch.chaput.ea:projected}.

Our results have been utilized by the second author to give an
explicit formula for the Schubert structure constants of the
equivariant quantum $K$-theory of projective space $\QK_T(\bP^n)$.
This formula establishes Griffeth-Ram positivity in this case
\cite{chung:fixme}.

\Theorem{sumcoef} was observed independently by Changzheng Li and
Leonardo Mihalcea, who also obtained proofs in some cases.  We thank
Li and Mihalcea for helpful discussions on this subject.

\section{Quantum $K$-theory}

In this section we briefly recall the definitions used in the
statements of our results, as well as the background required to prove
them.  A more detailed introduction to quantum $K$-theory can be found
in \cite[\S 2]{buch.chaput.ea:chevalley}.

Let $G$ be a semisimple complex linear algebraic group and fix a
maximal torus $T$, a Borel subgroup $B$, and a parabolic subgroup $P$
such that $T \subset B \subset P \subset G$.  Let $W = N_G(T)/T$ be
the Weyl group of $G$, let $W_P = N_P(T)/T$ be the Weyl group of $P$,
and let $W^P \subset W$ be the subset of minimal representatives for
the cosets in $W/W_P$.  Each element $w \in W^P$ defines the {\em Schubert
varieties\/} $X_w = \ov{Bw.P}$ and $X^w = \ov{B^-w.P}$ in the flag
variety $X = G/P$, where $B^- \subset G$ denotes the opposite Borel
subgroup defined by $B\cap B^- = T$.  We have $\dim(X_w) =
\codim(X^w,X) = \ell(w)$, where $\ell(w)$ is the length of $w$.  A
simple root $\ga$ of $G$ is called \emph{cominuscule} if the coefficient of
$\ga$ is one when the highest root is written as a linear combination
of simple roots.  The flag variety $X$ is \emph{cominuscule} if $W^P$
contains a single simple reflection $s_\ga$ defined by a cominuscule
simple root $\ga$.  We will assume this in what follows.  In
particular, we can identify $H_2(X;\Z) = \Z\, [X_{s_\ga}]$ with the
group of integers $\Z$.

Given a non-negative degree $d \in H_2(X;\Z)$ we let $\Mb_{0,n}(X,d)$
denote the Kontsevich moduli space of $n$-pointed stable maps to $X$
of degree $d$ and genus zero, see \cite{fulton.pandharipande:notes}.
This space is equipped with evaluation maps $\ev_i : \Mb_{0,n}(X,d)
\to X$ for $1 \leq i \leq n$.  Given any closed subvariety $Z \subset
X$, the {\em curve neighborhood\/} $\Gamma_d(Z) =
\ev_2(\ev_1^{-1}(Z))$ is the union of all connected rational curves of
degree $d$ in $X$ that meet $Z$.  It was proved in
\cite{buch.chaput.ea:finiteness} that, if $Z$ is a Schubert variety in
$X$, then so is $\Gamma_d(Z)$.  For $w \in W^P$ we let $w(-d) \in W^P$
denote the unique element for which $\Gamma_d(X^w) = X^{w(-d)}$.
Given two opposite Schubert varieties $X^u$ and $X_v$, the
corresponding {\em projected Gromov-Witten variety\/} is defined by
$\Gamma_d(X^u,X_v) = \ev_3(\ev_1^{-1}(X^u) \cap \ev_2^{-1}(X_v))$.
This is the union of all connected rational curves of degree $d$ that
meet both $X^u$ and $X_v$.  It was shown in
\cite{buch.chaput.ea:projected} that projected Gromov-Witten varieties
are also projected Richardson varieties as studied in
\cite{knutson.lam.ea:projections}, hence non-empty projected
Gromov-Witten varieties are unirational with rational singularities.
This generalizes the fact that any non-empty Richardson variety $X^u
\cap X_v$ is rational with rational singularities
\cite{richardson:intersections, ramanan.ramanathan:projective,
  ramanathan:schubert, brion:positivity}.  We let $\dist(X^u,X_v)$
denote the smallest degree $d$ for which $\Gamma_d(X^u,X_v) \neq
\emptyset$.

Let $K_T(X)$ denote the Grothendieck ring of $T$-equivariant algebraic
vector bundles on $X$.  Every $T$-stable closed subvariety $Z \subset
X$ defines a class $[\cO_Z] \in K_T(X)$.  If $Z$ is unirational with
rational singularities, then we have $\euler{X}([\cO_Z]) = 1 \in
\Gamma$, see \cite[Cor.~4.18(a)]{debarre:higher-dimensional}.  The
{\em Schubert classes\/} $\cO^w = [\cO_{X^w}]$ for $w \in W^P$ form a
basis of $K_T(X)$ as a module over the subring $\Gamma = K_T(\pt)$.
An alternative basis is provided by the $B$-stable Schubert classes
$\cO_w = [\cO_{X_w}]$.  Let $\Gamma\llbracket q \rrbracket$ denote the
ring of formal power series in a single variable $q$ with coefficients
in $\Gamma$.  The {\em equivariant quantum $K$-theory ring\/}
$\QK_T(X)$ is a $\Gamma\llbracket q \rrbracket$-algebra, which as a
module over $\Gamma\llbracket q \rrbracket$ is defined by $\QK_T(X) =
K_T(X) \otimes_\Gamma \Gamma\llbracket q \rrbracket$.  Givental
defined the product in $\QK_T(X)$ in terms of structure constants
obtained as polynomial expressions of Gromov-Witten invariants
\cite{givental:wdvv}.  In this paper we will use an alternative
construction from \cite{buch.chaput.ea:projected,
  buch.chaput.ea:chevalley}.  For $u,v \in W^P$, define a power series
in $\QK_T(X)$ by
\[
\cO^u \odot \cO_v = \sum_{d \geq 0}\, [\cO_{\Gamma_d(X^u,X_v)}]\, q^d \,.
\]
Let $\psi : \QK_T(X) \to \QK_T(X)$ be the unique $\Gamma\llbracket q
\rrbracket$-linear map defined by $\psi(\cO^w) = \cO^{w(-1)}$.  The
product in $\QK_T(X)$ is the unique $\Gamma\llbracket q
\rrbracket$-bilinear operator $\star$ defined by
\cite[Prop.~3.2]{buch.chaput.ea:chevalley}
\[
\cO^u \star \cO_v = (1 - q \psi)(\cO^u \odot \cO_v) \,.
\]

\section{Proof of Theorems 1, 2, and 3}

Let $\chi : \QK_T(X) \to \Gamma\llbracket q \rrbracket$ be the
$\Gamma\llbracket q \rrbracket$-linear extension of the Euler
characteristic map.  Since we have $\chi \psi = \chi$, \Theorem{dist}
follows from the calculation
\[
\begin{split}
  \chi(\cO^u \star \cO_v) \
  &= \ \chi(1-q\psi)(\cO^u \odot \cO_v)
  \ = \ (1-q)\chi(\cO^u \odot \cO_v) \\
  &= \ (1-q) \sum_{d \geq \dist(X^u,X_v)} q^d
  \ = \ q^{\dist(X^u,X_v)} \,.
\end{split}
\]

It follows from \cite[Thm.~1]{buch.chaput.ea:finiteness} that the
group $\QK_T^\poly(X) = K_T(X) \otimes_\Gamma \Gamma[q]$ of finite
power series is a subring of $\QK_T(X)$.  Let $\mu : \Gamma[q] \to
\Gamma$ be the ring homomorphism defined by $\mu(q) = 1$ and
$\mu(\alpha) = \alpha$ for $\alpha \in \Gamma$.  If we consider
$\Gamma$ as a module over $\Gamma[q]$ through this map, then the
composition $\wh\chi = \mu \chi : \QK_T^\poly(X) \to \Gamma$ is a
$\Gamma[q]$-linear map.  Since both of the sets $\{\cO^u : u \in
W^P\}$ and $\{\cO_v : v \in W^P\}$ are bases for $\QK_T^\poly(X)$ over
$\Gamma[q]$, it follows from the identity $\wh\chi(\cO^u \star \cO_v)
= 1 = \wh\chi(\cO^u) \cdot \wh\chi(\cO_v)$ that $\wh\chi$ is a
homomorphism of $\Gamma[q]$-algebras.  This proves \Theorem{ringhom}.
Finally, \Theorem{sumcoef} follows from applying $\wh\chi$ to both
sides of \eqn{strconst}.


\providecommand{\bysame}{\leavevmode\hbox to3em{\hrulefill}\thinspace}
\providecommand{\MR}{\relax\ifhmode\unskip\space\fi MR }
\providecommand{\MRhref}[2]{%
  \href{http://www.ams.org/mathscinet-getitem?mr=#1}{#2}
}
\providecommand{\href}[2]{#2}

\end{document}